\documentstyle[12pt,amstex,amscd,graphicx]{amsart}

\begin{document}

\newtheorem{theorem}{Theorem}[section]
\newtheorem{prop}[theorem]{Proposition}
\newtheorem{lemma}[theorem]{Lemma}
\newtheorem{cor}[theorem]{Corollary}
\newtheorem{defn}[theorem]{Definition}
\newtheorem{conj}[theorem]{Conjecture}
\newtheorem{claim}[theorem]{Claim}
\newtheorem{rem}{Remark}
\newcommand{\boundary}{\partial}
\newcommand{\C}{{\mathbb C}}
\newcommand{\integers}{{\mathbb Z}}
\newcommand{\natls}{{\mathbb N}}
\newcommand{\ratls}{{\mathbb Q}}
\newcommand{\reals}{{\mathbb R}}
\newcommand{\proj}{{\mathbb P}}
\newcommand{\lhp}{{\mathbb L}}
\newcommand{\tube}{{\mathbb T}}
\newcommand{\cusp}{{\mathbb P}}
\newcommand\AAA{{\mathcal A}}
\newcommand\BB{{\mathcal B}}
\newcommand\CC{{\mathcal C}}
\newcommand\DD{{\mathcal D}}
\newcommand\EE{{\mathcal E}}
\newcommand\FF{{\mathcal F}}
\newcommand\GG{{\mathcal G}}
\newcommand\HH{{\mathcal H}}
\newcommand\II{{\mathcal I}}
\newcommand\JJ{{\mathcal J}}
\newcommand\KK{{\mathcal K}}
\newcommand\LL{{\mathcal L}}
\newcommand\MM{{\mathcal M}}
\newcommand\NN{{\mathcal N}}
\newcommand\OO{{\mathcal O}}
\newcommand\PP{{\mathcal P}}
\newcommand\QQ{{\mathcal Q}}
\newcommand\RR{{\mathcal R}}
\newcommand\SSS{{\mathcal S}}
\newcommand\TT{{\mathcal T}}
\newcommand\UU{{\mathcal U}}
\newcommand\VV{{\mathcal V}}
\newcommand\WW{{\mathcal W}}
\newcommand\XX{{\mathcal X}}
\newcommand\YY{{\mathcal Y}}
\newcommand\ZZ{{\mathcal Z}}
\newcommand\CH{{\CC\HH}}
\newcommand\TC{{\TT\CC}}
\newcommand\EXH{{ \EE (X, \HH )}}
\newcommand\GXH{{ \GG (X, \HH )}}
\newcommand\GYH{{ \GG (Y, \HH )}}
\newcommand\PEX{{\PP\EE  (X, \HH , \GG , \LL )}}
\newcommand\MF{{\MM\FF}}
\newcommand\PMF{{\PP\kern-2pt\MM\FF}}
\newcommand\ML{{\MM\LL}}
\newcommand\PML{{\PP\kern-2pt\MM\LL}}
\newcommand\GL{{\GG\LL}}
\newcommand\Pol{{\mathcal P}}
\newcommand\half{{\textstyle{\frac12}}}
\newcommand\Half{{\frac12}}
\newcommand\Mod{\operatorname{Mod}}
\newcommand\Area{\operatorname{Area}}
\newcommand\ep{\epsilon}
\newcommand\hhat{\widehat}
\newcommand\Proj{{\mathbf P}}
\newcommand\U{{\mathbf U}}
 \newcommand\Hyp{{\mathbf H}}
\newcommand\D{{\mathbf D}}
\newcommand\Z{{\mathbb Z}}
\newcommand\R{{\mathbb R}}
\newcommand\Q{{\mathbb Q}}
\newcommand\E{{\mathbb E}}
\newcommand\til{\widetilde}
\newcommand\length{\operatorname{length}}
\newcommand\tr{\operatorname{tr}}
\newcommand\gesim{\succ}
\newcommand\lesim{\prec}
\newcommand\simle{\lesim}
\newcommand\simge{\gesim}
\newcommand{\simmult}{\asymp}
\newcommand{\simadd}{\mathrel{\overset{\text{\tiny $+$}}{\sim}}}
\newcommand{\ssm}{\setminus}
\newcommand{\diam}{\operatorname{diam}}
\newcommand{\pair}[1]{\langle #1\rangle}
\newcommand{\T}{{\mathbf T}}
\newcommand{\inj}{\operatorname{inj}}
\newcommand{\pleat}{\operatorname{\mathbf{pleat}}}
\newcommand{\short}{\operatorname{\mathbf{short}}}
\newcommand{\vertices}{\operatorname{vert}}
\newcommand{\collar}{\operatorname{\mathbf{collar}}}
\newcommand{\bcollar}{\operatorname{\overline{\mathbf{collar}}}}
\newcommand{\I}{{\mathbf I}}
\newcommand{\tprec}{\prec_t}
\newcommand{\fprec}{\prec_f}
\newcommand{\bprec}{\prec_b}
\newcommand{\pprec}{\prec_p}
\newcommand{\ppreceq}{\preceq_p}
\newcommand{\sprec}{\prec_s}
\newcommand{\cpreceq}{\preceq_c}
\newcommand{\cprec}{\prec_c}
\newcommand{\topprec}{\prec_{\rm top}}
\newcommand{\Topprec}{\prec_{\rm TOP}}
\newcommand{\fsub}{\mathrel{\scriptstyle\searrow}}
\newcommand{\bsub}{\mathrel{\scriptstyle\swarrow}}
\newcommand{\fsubd}{\mathrel{{\scriptstyle\searrow}\kern-1ex^d\kern0.5ex}}
\newcommand{\bsubd}{\mathrel{{\scriptstyle\swarrow}\kern-1.6ex^d\kern0.8ex}}
\newcommand{\fsubeq}{\mathrel{\raise-.7ex\hbox{$\overset{\searrow}{=}$}}}
\newcommand{\bsubeq}{\mathrel{\raise-.7ex\hbox{$\overset{\swarrow}{=}$}}}
\newcommand{\tw}{\operatorname{tw}}
\newcommand{\base}{\operatorname{base}}
\newcommand{\trans}{\operatorname{trans}}
\newcommand{\rest}{|_}
\newcommand{\bbar}{\overline}
\newcommand{\UML}{\operatorname{\UU\MM\LL}}
\newcommand{\EL}{\mathcal{EL}}
\newcommand{\tsum}{\sideset{}{'}\sum}
\newcommand{\tsh}[1]{\left\{\kern-.9ex\left\{#1\right\}\kern-.9ex\right\}}
\newcommand{\Tsh}[2]{\tsh{#2}_{#1}}
\newcommand{\qeq}{\mathrel{\approx}}
\newcommand{\Qeq}[1]{\mathrel{\approx_{#1}}}
\newcommand{\qle}{\lesssim}
\newcommand{\Qle}[1]{\mathrel{\lesssim_{#1}}}
\newcommand{\simp}{\operatorname{simp}}
\newcommand{\vsucc}{\operatorname{succ}}
\newcommand{\vpred}{\operatorname{pred}}
\newcommand\fhalf[1]{\overrightarrow {#1}}
\newcommand\bhalf[1]{\overleftarrow {#1}}
\newcommand\sleft{_{\text{left}}}
\newcommand\sright{_{\text{right}}}
\newcommand\sbtop{_{\text{top}}}
\newcommand\sbot{_{\text{bot}}}
\newcommand\sll{_{\mathbf l}}
\newcommand\srr{_{\mathbf r}}
\newcommand\geod{\operatorname{\mathbf g}}
\newcommand\mtorus[1]{\boundary U(#1)}
\newcommand\A{\mathbf A}
\newcommand\Aleft[1]{\A\sleft(#1)}
\newcommand\Aright[1]{\A\sright(#1)}
\newcommand\Atop[1]{\A\sbtop(#1)}
\newcommand\Abot[1]{\A\sbot(#1)}
\newcommand\boundvert{{\boundary_{||}}}
\newcommand\storus[1]{U(#1)}
\newcommand\Momega{\omega_M}
\newcommand\nomega{\omega_\nu}
\newcommand\twist{\operatorname{tw}}
\newcommand\modl{M_\nu}
\newcommand\MT{{\mathbb T}}
\newcommand\Teich{{\mathcal T}}
\renewcommand{\Re}{\operatorname{Re}}
\renewcommand{\Im}{\operatorname{Im}}

\title{Relative Hyperbolicity, Trees of Spaces and
 Cannon-Thurston Maps}

\author{Mahan Mj}
\address{RKM Vivekananda University, Belur Math, WB-711 202, India}

\author{Abhijit Pal}
\address{Stat-Math Unit, Indian Statistical Institute, 203 B.T.Road, Kolkata 700108}

\thanks{This paper is part of AP's PhD thesis  written under the supervision of MM}

\date{}

\begin{abstract}
We prove the existence of continuous boundary extensions (Cannon-Thurston maps) for the inclusion of a vertex space into a tree of (strongly) relatively
hyperbolic spaces satisfying the qi-embedded condition. This
implies the same  result for
inclusion of vertex (or edge)
subgroups in finite graphs of (strongly) relatively
hyperbolic groups. This generalizes a result of Bowditch for punctured surfaces in 3 manifolds and a result of Mitra for
trees of hyperbolic metric spaces.
\end{abstract}

\maketitle

\begin{center}
AMS subject classification =   20F32(Primary), 57M50(Secondary)
\end{center}

\tableofcontents

\section{Introduction}

For a closed hyperbolic 3-manifold $M$, fibering over the circle with
fiber $F$, let $i: \widetilde{F} \rightarrow \widetilde{M}$ denote the inclusion of universal
covers. In \cite{CT} (now published as \cite{CTpub})
 Cannon and Thurston
show that $i$
extends to a continuous map  $\hat{i} : {\Bbb{D}}^2 \rightarrow {\Bbb{D}}^3$ where
${{\Bbb{D}}^2}={\Bbb{H}}^2\cup{\Bbb{S}}^1_\infty$ and
${{\Bbb{D}}^3}={\Bbb{H}}^3\cup{\Bbb{S}}^2_\infty$
denote the standard compactifications. In \cite{minsky-jams}, Minsky generalized Cannon and Thurston's result to
bounded geometry surface Kleinian
groups {\em without} parabolics

In \cite{mitra-trees}, one of us extended Cannon-Thurston's and Minsky's  result to trees of hyperbolic metric spaces satisfying a natural qi-embedded  condition.
In the process, an alternate proof of Cannon-Thurston's original result was found.

Bowditch
\cite{bowditch-ct} \cite{bowditch-stacks} made use of
 some of the ideas of \cite{mitra-trees} amongst other things and proved
 the Cannon-Thurston
property for bounded geometry surface Kleinian
groups with parabolics. (It is worth remarking parenthetically that this generalization from the case without punctures to that with punctures required essentially new ideas and a fair bit of time.)

In \cite{brahma-bddgeo} one of us gave a different proof of Bowditch's result.

The appropriate framework
 for synthesizing and generalizing the above results is that of trees of (strong) relatively hyperbolic metric spaces. A combination theorem was described by Mj and
Reeves in \cite{mahan-reeves}. The notion of partial electrocution introduced there will be used essentially here.
Relatively hyperbolic spaces (Gromov \cite{gromov-hypgps}, Farb \cite{farb-relhyp}, Bowditch \cite{bowditch-relhyp}, etc.)
generalize fundamental groups of finite volume manifolds of pinched negative curvature. We shall implicitly use the fact, due to Bowditch \cite{bowditch-relhyp}, that the (strong) relative hyperbolic
boundary of a space is well-defined. Unless otherwise mentioned,
relative hyperbolicity will mean strong relative hyperbolicity.
Our main Theorem is: \\

\medskip

\noindent {\bf Theorem  \ref{main thm}:}
 Let $P\colon X\to T$ be a tree  of relatively hyperbolic spaces satisfying
the quasi-isometrically (qi) embedded condition. Further suppose that inclusion of edge-spaces into vertex spaces is strictly
type-preserving, and that the induced tree of coned-off spaces continues to satisfy the qi-embedded condition.
If $X$  is {\em
  strongly hyperbolic} relative to the family $\CC$ of maximal
  cone-subtrees of horosphere-like sets,
then a Cannon-Thurston map exists for the proper
embedding $i\colon X_{v}\to X$, where
$v$ is a vertex of $T$ and
$({X_v},d_{X_v})$ is the relatively hyperbolic metric space corresponding to $v$.

\medskip

A special case of Theorem  \ref{main thm} is that of punctured surface Kleinian groups, where we obtain
the following due to Bowditch \cite{bowditch-ct}.

\medskip

 \noindent {\bf Theorem  \ref{bowditch}:}({\bf Bowditch \cite{bowditch-ct}})
Let $M$ be a 3-manifold corresponding to a Kleinian surface group without accidental parabolics. Further, suppose that
$M$ has bounded geometry. If $S$ denotes the corresponding finite volume hyperbolic surface with some hyperbolic structure, then the inclusion $i: \til{S} \rightarrow
\til{M}$ extends continuously to the boundary, i.e. has a
Cannon-Thurston map.

\medskip

In fact, we obtain the following more general theorem due to the first author \cite{brahma-pared}
for $M$ any bounded
geometry hyperbolic 3 manifold with core incompressible
away from cusps (not necessarily a surface group), and {\it no accidental parabolics}:

\medskip

\noindent {\bf Theorem  \ref{pared}:} \cite{brahma-pared}
Let $M$ be a hyperbolic
3-manifold corresponding to a Kleinian group of bounded geometry without
accidental parabolics. Further suppose that the compact core
of $M$ is incompressible away from cusps. If $N$ denotes
a geometrically finite hyperbolic
3 manifold with some hyperbolic structure identified with
the convex core of $M$, then
the inclusion $i: \til{N} \rightarrow
\til{M}$ extends continuously to the boundary, i.e. has a
Cannon-Thurston map.

\medskip

\noindent {\bf Acknowledgements:} We would like to thank the referee for detailed and helpful comments.
The research of the first author is partly supported by a Department
of Science and Technology Research Project grant.

\subsection{Hyperbolicity and Nearest Point Projections}
We assume that the reader is familiar with the basic notions about hyperbolic metric spaces in the sense of Gromov \cite{gromov-hypgps}. (See \cite{GhH} for instance).
For a hyperbolic metric space $X$, the Gromov bordification will be denoted by $\bbar{X}$.
Here, we give some basic facts about hyperbolic metric spaces.

\medskip

\begin{defn}
 A path $\gamma \colon I\to X$ in a path metric space $X$ is said to be a {\bf  $K$-quasigeodesic} if we have
\[
 L(\beta)\leq KL(A)+K
\]
for any subsegment $\beta = \gamma |[a,b]$ and any rectifiable path $A\colon [a,b]\to X$ with the same end points. Here $L(\sigma)$ denotes the length of a path $\sigma$.
\end{defn}

For points $x,y$ in a geodesic metric space $X$, let $[x,y]$ denote a geodesic in $X$ joining $x$ and $y$.  The following is an easy consequence of $\delta$-hyperbolicity.

\begin{lemma}\label{gen lem} (Lemma $3.1$ of \cite{mitra-trees}) \\
Given $\delta > 0$,  there exist $D, C_1$ such that if $x, y$ are points
of a $\delta -$ hyperbolic metric space $(X, d)$, $\lambda$ is a hyperbolic geodesic in $X$, and
  $\pi_{\lambda}$ is a nearest point projection
of $X$ onto $\lambda$ with $d(\pi_{\lambda}(x), \pi_{\lambda}(y))\geq D$,   then $[x,\pi_{\lambda}(x)]\cup [\pi_{\lambda}(x),\pi_{\lambda}(y)]\cup [\pi_{\lambda}(y), y] $ lies in a $C_1-$neighborhood of any geodesic joining $x, y$.
\end{lemma}

The next Lemma states that nearest point projections in a $\delta$-hyperbolic metric spaces  do not increases distance much.
\begin{lemma}\label{proj lem}(See Lemma 3.2 of \cite{mitra-trees})
For a $\delta$-hyperbolic metric space $(Y, d)$, let $\pi_\lambda$ be the nearest point projection onto
the geodesic segment $\lambda$.  There exists $P_1>0$ (depending only on $\delta$) such that
$d(\pi_\lambda(x), \pi_\lambda(y))\leq P_1d(x, y)+ P_1$ for all $x, y\in Y$.
\end{lemma}

The following Lemma says that nearest point projections and quasiisometries in  hyperbolic metric spaces  \textquoteleft almost commute\textquoteright.

\begin{lemma}(See Lemma 3.5 of \cite{mitra-trees})\label{q.i. lem}
Suppose $(Y_1, d_1)$ and $(Y_2, d_2)$ are $\delta$-hyperbolic metric spaces.
Let $\mu_1$ be some geodesic segment in $Y_1$ joining $a, b$ and let $p\in Y_1$.
Let $\phi$ be a $(K, {\epsilon})$ - quasiisometry from $Y_1$ to $Y_2$.
Let $\mu_2$ be a geodesic segment in $Y_2$ joining ${\phi}(a)$ to ${\phi}(b)$.
Then $d_{Y_2}(\pi_{\mu_2}(\phi (p)), \phi (\pi_{\mu_1}(p)))\leq{P_2}$ for some constant $P_2$ dependent only on
$K,  \epsilon $ and $\delta$.
\end{lemma}

{\bf Remark:} Due to stability of quasigeodesics, the Lemmas \ref{proj lem} and
\ref{q.i. lem} are also true if geodesics are replaced by quasigeodesics and  nearest point projections are taken onto  quasigeodesics.
For  ease of exposition, we will use the
same constants as above when geodesics are replaced by quasigeodesics.

\subsection{Relative Hyperbolicity}

Let $(X,d)$ be a path metric space. A collection of closed
 subsets $\HH = \{ H_\alpha\}$ of $X$ will be said to be {\bf uniformly
 separated} if there exists $\epsilon > 0$ such that
$d(H_1, H_2) \geq \epsilon$ for all distinct $H_1, H_2 \in \HH$.

\begin{defn} (Farb \cite{farb-relhyp})
The {\bf electric space} (or coned-off space) $\EE {(X, \HH )}$
corresponding to the
pair $(X,\HH )$ is a metric space which consists of $X$ and a
collection of vertices  $v_\alpha$ (one for each $H_\alpha \in \HH$)
such that each point of $H_\alpha$ is joined to (coned off at)
$v_\alpha$ by an edge of length $\half$. The sets $H_\alpha$ shall be
referred to as  horosphere-like sets and the vertices $v_{\alpha}$ as cone-points. \\
$X$ is said to be {\bf weakly hyperbolic} relative to the collection $\HH$ if $\EE {(X, \HH )}$
is a hyperbolic metric space.
\label{el-space}
\end{defn}

\begin{defn}
$\bullet$ A path $\gamma$ in $\EXH$ is said to be an electric geodesic (resp. electric $K$-quasigeodesic) if it is a geodesic (resp. $K$-quasigeodesic) in $\EXH$.\\
$\bullet$ $\gamma$ is said to be an  electric
$K$-quasigeodesic in (the electric space) $\EE {(X, \HH )}$
 without  backtracking  if
 $\gamma$ is an electric $K$-quasigeodesic in $\EE {(X, \HH )}$ and
 $\gamma$ does not return to  any {\em horosphere-like set} $H_\alpha$ after leaving it.
\end{defn}

Let $i: X \rightarrow \EE {(X, \HH )}$ denotes the natural inclusion of spaces. Then
for a path $\gamma \subset X$, the path $i ( \gamma )$ lies in $\EE {(X, \HH )}$. Replacing maximal subsegments
$[a,b]$ of $i (\gamma )$  lying in a particular $H_\alpha$ by a path that goes from
$a$ to $v_\alpha$ and then from $v_\alpha$ to $b$, and repeating this for every $H_\alpha$ that
$i ( \gamma )$ meets we obtain a new path $\hat{\gamma}$. If $\hat{\gamma}$ is an electric geodesic (resp. $P$-quasigeodesic), $\gamma$ is called a {\em relative geodesic} (resp.
{\em relative $P$-quasigeodesic}). We shall usually be concerned with the case that $\gamma$ is an ambient geodesic/quasigeodesic without backtracking.

\begin{defn}
Relative  $P$-quasigeodesics in
$(X,\HH )$ are said to satisfy {\bf bounded region penetration} if for any two
    relative  $P$-quasigeodesics without backtracking
$\beta$, $\gamma$,
   joining $x, y \in X$
 there exists $B = B(P )$ such that \\
{\bf Similar Intersection Patterns 1:}  if
  precisely one of $\{ \beta , \gamma \}$ meets
 a  horosphere-like set $H_\alpha$,
then the length (measured in the intrinsic path-metric
  on  $H_\alpha$ ) from the first (entry) point
  to the last
  (exit) point (of the relevant path) is at most $B$. \\
 {\bf Similar Intersection Patterns 2:}  if
 both $\{ \beta , \gamma \}$ meet some  $H_\alpha $
 then the length (measured in the intrinsic path-metric
  on  $H_\alpha$ ) from the entry point of
 $\beta$ to that of $\gamma$ is at most $B$; similarly for exit points. \\
\end{defn}

Replacing `$P$-quasigeodesic' by `geodesic' in the above definition, we obtain the notion of
relative  geodesics in
$(X,\HH )$   satisfying  bounded region penetration.

Families of paths which enjoy the above properties shall be said to have similar intersection patterns with horospheres.

\begin{defn} (Farb \cite{farb-relhyp} ) $X$ is said to be hyperbolic relative to the uniformly separated collection $\HH$ if \\
1) $X$ is weakly hyperbolic relative to $\HH$ \\
2) For all $P \geq 1$, relative $P$-quasigeodesics without backtracking satisfy the bounded penetration property \\
\end{defn}

\noindent {\underline{Gromov's definition of  relative hyperbolicity \cite{gromov-hypgps} }:}\\
 \begin{defn} ({\bf Gromov})
For any geodesic metric space
$(H,d)$, the {\em hyperbolic cone} (analog of a horoball)
$H^h$ is the metric space
$H\times [0,\infty) = H^h$ equipped with the
path metric $d_h$ obtained from two pieces of
 data \\
1) $d_{h,t}((x,t),(y,t)) = 2^{-t}d_H(x,y)$, where $d_{h,t}$ is the induced path
metric on $H\times \{t\}$.  Paths joining
$(x,t),(y,t)$ and lying on  $H\times \{t\}$
are called {\em horizontal paths}. \\
2) $d_h((x,t),(x,s))=\vert t-s \vert$ for all $x\in H$ and for all $t,s\in [0,\infty)$, and the corresponding paths are called
{\em vertical paths}. \\
3)  for all $x,y \in H^h$,  $d_h(x,y)$ is the path metric induced by the collection of horizontal and vertical paths. \\
\end{defn}

\begin{defn}
Let $X$ be a geodesic metric space and $\HH$ be a collection of mutually disjoint uniformly separated subsets of $X$.
$X$ is said to be  hyperbolic relative to $\HH$ in the sense of Gromov, if the quotient space $\GG (X, \HH)$,  obtained by attaching the hyperbolic cones
$ H^h$ to $H \in \HH$  by identifying $(z,0)$ with $z$
for all $H\in \HH$ and $z \in H$,
 is a complete hyperbolic metric space. The collection $\{ H^h : H \in \HH \}$ is denoted
as ${\HH}^h$. The induced path metric is denoted as $d_h$.
\end{defn}

Here $H\in \HH$ are thought of as horosphere-like sets and $H\times [0,\infty)$ as horoballs.
For a strong relatively hyperbolic metric space $(X,\HH)$, the space $\EE (\GXH , \HH^h)$ obtained by coning off $H\times [0,\infty)$ for all $H\in \HH$ is basically the same as $\EE (X, {\HH})$. The latter embeds isometrically
into the former and every point of $\EE (\GXH , \HH^h)$ lies within a distance one of the isometric image
of $\EE (X, {\HH})$. Thus we can (and shall) pass interchangeably between $\EXH$
and $\EE (\GXH , \HH^h)$.

\begin{theorem} (Bowditch \cite{bowditch-relhyp})\label{bow-rel}
The following are equivalent: \\
1) $X$ is  hyperbolic relative to the
collection $\HH$ of uniformly separated subsets of $X$  \\
2) $X$ is  hyperbolic relative to the
collection $\HH$ of uniformly separated subsets of $X$ in the sense of Gromov \\
3) $\GXH$ is  hyperbolic relative to the
collection $\HH^h$ \\
\end{theorem}

We collect together certain facts about the electric metric that Farb
proves in \cite{farb-relhyp}. These are proved in the context of Hadamard manifolds of pinched negative curvature, but the proofs go through in our context. $N^h_R(Z)$ will denote the
$R$-neighborhood about the subset $Z$ in $(\GXH , d_h)$.
 $N_R^e(Z)$ will denote the
$R$-neighborhood about the subset $Z$ in the electric metric $(\EXH , d_e)$. Geodesics
in $(\GXH , d_h)$ will be termed $h$-geodesics. Similarly for quasigeodesics.

\begin{lemma} (See Lemma 4.5 and Proposition 4.6 of \cite{farb-relhyp})
\begin{enumerate}
\item {\it Electric quasigeodesics electrically track hyperbolic
  geodesics:} Given $P > 0$, there exists $K > 0$ with the following
  property: \\
Let $\beta$ be any electric $P$-quasigeodesic without backtracking
from $x$ to
  $y$ in $\EXH$, and let $\gamma$ be an $h$-geodesic from $x$ to $y$ in $\GXH$.
Then $\beta \subset N_K^e ( \gamma )$.
\item {\it Quasiconvexity:} There exists $K$ such that each $H^h$
 is uniformly $K$-quasiconvex in $\GXH$.
\item electric geodesics in $\EXH$ and
 relative geodesics in $X$
joining the same pair of points in $X$ have similar intersection patterns with $H$
for all $H \in \HH$, i.e. they track each other off horosphere-like  sets.
\item electric geodesics in $\EXH$ (after identification with $\EE ( \GXH , \HH^h )$) and
 $h$-geodesics in $\GXH$   joining the same pair
of points in  $\GXH$ have similar intersection patterns with $H^h$
for all $H^h \in \HH^h$, i.e. they track each other off horoball-like  sets.
\end{enumerate}
\label{farb1}
\end{lemma}
\begin{defn} Let $X$ be hyperbolic relative to $\HH$.
We start with an electric quasi-geodesic $\hat \lambda$ in the electric space $\EXH$ without backtracking.   For any  $H\in \HH$ penetrated by $\hat \lambda$,  let $x_H$ and $y_H$ be the first entry
point  and the last exit point of $\hat \lambda$.   We join $x_H$ and $y_H$ by a hyperbolic geodesic segment in $H^h$ (identifying $\EXH$ with  $\EE ( \GXH, \HH^h )$).   This results in a path $\lambda$
  in $\GXH$.   The path $\lambda$ will be called an {\bf electro-ambient quasigeodesic}.
\end{defn}

\begin{lemma}\label{ea}
An electro-ambient quasigeodesic is a quasigeodesic in  $\GXH$.
\end {lemma}

\begin{center}
\bf{Electric Projections}
\end{center}

Let $Y$ be a space hyperbolic relative to the collection $\HH_Y$, $\widehat Y= \EE(Y,\HH_Y)$, $Y^h=\GG(Y,\HH_Y)$, $\delta$ be the common hyperbolicity constant and  $i\colon Y^h\to \widehat{Y^h}=\EE(Y^h,\HH^h_Y)$ denote the inclusion. Here, as pointed
out earlier, we identify $\widehat{Y^h}$ with $\widehat Y$.
Suppose $\hat\mu$ is a geodesic in $\widehat Y$, $\mu$ is an electro-ambient representative of the geodesic $\hat \mu$ in $Y^h$
and $\pi_\mu$ is a nearest point projection from $Y^h$ onto $\mu$.

\begin{defn}(Electric Projection)\label{proj}
Let $y\in \widehat Y$ and $\hat\mu$ be a geodesic in $\widehat Y$. 
Define $\hat \pi_{\hat \mu} : \widehat{Y} \rightarrow
\mu$ as follows.\\For $y\in Y$,  $\hat \pi_{\hat \mu}(y)=i(\pi_\mu(y))$. \\If $y$ is a cone point over a horosphere-like set $H\in \HH_Y$, choose $z\in H$ and define $\hat \pi_{\hat \mu}(y)=i(\pi_\mu(z))$.\\
$\hat \pi_{\hat \mu}$ will be called an Electric Projection.
\end{defn}

The next lemma, shows that $\hat \pi_{\hat \mu}$ is well-defined up to a bounded amount of discrepancy with respect to the metric $d_{\widehat Y}$. Let $D, C_1$ be as in Lemma \ref{gen lem}.

\begin{lemma}\label{bdd disc}
Let $Y$ be  hyperbolic relative to the collection $\HH_Y$. There exists a constant $P_3>0$ depending only upon $\delta, D, C_1$ such that for any $H\in \HH_Y$ and $z,z'\in H$, if $\hat\mu$ be a geodesic in $\widehat Y$ then $d_{\widehat Y}(i(\pi_\mu(z)),i(\pi_\mu(z')))\leq P_3$.
\end{lemma}
{\bf Proof.} If $d_{Y^h}(\pi_\mu(z),\pi_\mu(z'))\leq D$, then $d_{\widehat Y}(i(\pi_\mu(z)),i(\pi_\mu(z')))\leq D$.\\
Next suppose $d_{Y^h}(\pi_\mu(z),\pi_\mu(z'))> D$. Let $a = \pi_{\lambda}(z)$ and $b = \pi_{\lambda}(z')$, then $[z,a]\cup [a,b]\cup [b,z']=\beta (say)$ is a  quasi-geodesic in $Y^h$. Hence $\hat\beta$ is an electric $(L,\epsilon_1)$ quasi-geodesic in $\widehat Y$.
Thus
\begin{eqnarray*}
d_{\widehat Y}(i(\pi_\mu(z)),i(\pi_\mu(z')))
&\leq&l_{\widehat Y}(\hat\beta)\\
&\leq& L d_{\widehat Y}(z,z') + \epsilon_1\\
&\leq& L + \epsilon_1 = D_1 (say),\mbox{ since  }d_{\widehat Y}(z,z')\leq 1.
\end{eqnarray*}
Let $P_3=max\{D_1,D'\}$, then we have the required result. $\Box$\\
\medskip

Further, if $x,y\in\widehat Y$ and $d_{\widehat Y}(x,y)\leq 1$ then similarly  we can prove that there exists constant $R>0$ (depending only upon $\delta, D, C_1$) such that $d_{\widehat Y}(\hat \pi _{\hat \lambda}(x),\hat \pi _{\hat \lambda}(y))\leq R$.
Thus we have the following lemma:
\begin{lemma}\label{ret lem}
Let $Y$ be  hyperbolic relative to the collection $\HH_Y$. For the electric space $\widehat Y$ there exists $P_4>0$ (depending only upon the hyperbolic constant of $\widehat Y$) such that for all $x,y\in \widehat Y$ and electric geodesic segment $\hat \lambda$,
$d_{\widehat Y}(\hat \pi _{\hat \lambda}(x),\hat \pi _{\hat \lambda}(y))\leq P_4 d_{\widehat Y}(x,y)+P_4$.
\end{lemma}
{\bf Note:} Electric projection may not be a nearest point projection from an electric space onto an electric geodesic but in analogy with Lemma \ref{proj lem}, the above lemma says that electric projections do not increases distance much.

\medskip

As an easy
 consequence of Lemma \ref{q.i. lem}, we have the following Lemma which says that electric projections and \textquoteleft{strictly type-preserving\textquoteright}  quasiisometries  \textquoteleft{almost commute\textquoteright} in electric spaces.
\begin{lemma}\label{e.qi}
Suppose $Y_1$ and $Y_2$ are two metric spaces hyperbolic relative to the collections $\HH_{Y_1}$ and $\HH_{Y_2}$ respectively and $\phi\colon Y_1\to Y_2$ is a $(K,\epsilon)$- strictly type-preserving quasiisometry, i.e., for all $H_{Y_1}\in\HH_{Y_1},H_{Y_2}\in\HH_{Y_2}$, $\phi(H_{Y_1})\in\HH_{Y_2}$ and $\phi^{-1}(H_{Y_2})\in\HH_{Y_1}$. Let $\hat\mu_1$ be a quasigeodesic in $\widehat Y_1$ joining $a, b$ and $\hat\phi$ be the induced  quasiisometry from $\widehat Y_1$ to $\widehat Y_2$.
Let $\hat\mu_2$ be a quasigeodesic in $\widehat Y_2$ joining ${\hat\phi}(a)$ to ${\hat\phi}(b)$. If $p\in\widehat Y_1$ then $d_{\widehat Y_2}(\hat\pi_{\hat\mu_2}(\hat\phi (p)), \hat\phi (\hat\pi_{\hat\mu_1}(p)))\leq{P_5}$, for some constant $P_5$ dependent only on $\delta,K,\epsilon$.
\end{lemma}

\subsection{Partial Electrocution}
In this subsection, we summarize some material   from
\cite{mahan-reeves} by Mj and Reeves.

\begin{defn}
Let $(X, \HH , \GG , \LL )$ be an ordered quadruple such that the
following holds:\\
\begin{enumerate}
\item $X$ is a geodesic metric space. $\HH$ is a collection of subsets $H_\alpha$ of $X$.
 $X$ is  hyperbolic relative to $\HH$.
\item There exists $\delta > 0$ such that $\LL$ is a collection of $\delta$-hyperbolic metric spaces
$L_\alpha$
and $\GG$ is a collection of (uniformly) coarse Lipschitz maps
$g_\alpha : H_\alpha \rightarrow L_\alpha$. Note that the indexing set for $H_\alpha, L_\alpha, g_\alpha$ is common.
\end{enumerate}
 The {\bf partially electrocuted space} or
{\em partially coned off space}  $\PEX $
corresponding to $(X, \HH , \GG , \LL)$
is
obtained from $X$ by gluing in the (metric)
mapping cylinders for the maps
 $g_\alpha : H_\alpha \rightarrow L_\alpha$. The metric on $\PEX$ is denoted by
$d_{pel}$.
\label{pex}
\end{defn}

In Farb's construction \cite{farb-relhyp}, each
 $L_\alpha$ is just a single point and $g_\alpha$ a constant map and we recover the
definition of an electric space in the sense of Farb from the above definition. We thus
think of Farb's definition as that of a `completely electrocuted space'.

 The metric,  geodesics and quasigeodesics
in the partially electrocuted space will be referred to as the
partially electrocuted metric $d_{pel}$, and partially
electrocuted geodesics and quasigeodesics respectively.

The generalization of the technical core of \cite{farb-relhyp} to our context
requires reworking the arguments in \cite{farb-relhyp} to our context. The proofs of Lemmas \ref{pel}
and \ref{pel-track} below will be given together after stating Lemma \ref{pel-track}.

\begin{lemma}
For $(X, \HH , \GG , \LL )$  an ordered quadruple as in Definition
\ref{pex} above,
$(\PEX ,d_{pel})$ is a hyperbolic metric space and the sets $L_\alpha$
are uniformly quasiconvex.
\label{pel}
\end{lemma}

\begin{lemma}
Let $(X, \HH , \GG , \LL )$  be an ordered quadruple as in Definition
\ref{pex} above.
Given $K, \epsilon \geq 0$, there exists $C > 0$ such that the following
holds: \\
Let $\gamma_{pel}$ and $\gamma$ denote respectively a $(K, \epsilon )$
partially electrocuted quasigeodesic in $(\PEX,d_{pel})$ and a
$(K, \epsilon )$ h-quasigeodesic in $(\GXH ,d_h)$ joining $a, b$. Then $\gamma \setminus
\bigcup_{H_\alpha\in\HH} H_\alpha$
lies in a  $C$-neighborhood of (any representative of)
$\gamma_{pel}$ in $(X,d)$. Further, outside of  a $C$-neighborhood of the horoballs
that $\gamma$ meets, $\gamma$ and $\gamma_{pel}$ track each other, i.e. lie in a $C$-neighborhood
of each other.
\label{pel-track}
\end{lemma}

\noindent {\bf Proofs of Lemmas \ref{pel}
and \ref{pel-track} :}
 By Theorem \ref{bow-rel}, a relatively hyperbolic space $(X, \HH )$
(in Farb's sense) can
be relatively hyperbolized (in Gromov's sense) to $\GXH$. Note that $(\PEX,d_{pel})$ is (strongly)
 hyperbolic relative to
the sets $\{ L_\alpha \}$. In fact the space obtained by electrocuting the
sets $L_\alpha$ in $(\PEX ,d_{pel})$ is just the space $(\EXH,d_e)$. We thus identify the three spaces
$\EE (X, \HH )$, $\EE (\GXH , {\HH}^h)$, $\EE ( \PEX , \LL )$.

Hence by part (4) of Lemma \ref{farb1}  geodesics in $(\EE (X, \HH ), d_e)$,
$(\GXH , d_h)$ and $ ( \PEX , d_{pel} )$
track geodesics in $(X, d)$ off horosphere-like sets.

By part (4) of Lemma \ref{farb1} applied to $\EE (\GXH , {\HH}^h)$ again, each element of ${\HH}^h$ is uniformly
quasiconvex in $\GXH $. Now, let $a, b \in (\GXH \setminus \bigcup_{H^h \in {\HH}^h}H^h)$ and let $Z^h_{a,b}$ denote
the union of the geodesic $\lambda_{a,b}$
in $\GXH$ joining $a, b$ along with the elements of $\HH^h$ it intersects nontrivially. Also let $Z_{a,b}
= Z^h_{a,b} \cap X$ be the union of the part of $\lambda_{a,b}$ lying in $X$ along with the {\it horosphere-like
sets} (elements of $\HH$) it intersects.

It is easy to check that $Z^h_{a,b}$ is uniformly
quasiconvex in $\GXH $. (For instance, if $c, d$ are two points in two different horoballs in $Z^h_{a,b}$,
join them to $\lambda_{a,b}$ by the shortest geodesic `perpendicular'
lying inside these horoballs and adjoin the segment
of $\lambda_{a,b}$ between the feet of the `perpendiculars'. The union of these three
segments is a uniform quasigeodesic.) Hence  
 the nearest-point projection of $\GXH$ onto   $Z^h_{a,b}$ decreases  $d_h$
distances by an exponential factor. Restricting this projection to
$X$  it follows that $Z_{a,b}$ is quasiconvex in $X$ and the nearest point projection of $X$ onto
$Z_{a,b}$ decreases the $d$
distance by an exponential factor. 

By the identification of $\EE (X, \HH )$ and $ \EE ( \PEX , \LL )$ mentioned above, and the (strong)
hyperbolicity of 
$(\PEX,d_{pel})$ 
  relative to
the sets $\{ L_\alpha \}$, it follows that $Z_{a,b} \subset \PEX$ (thought of as a subset of $\PEX$ under the natural
inclusion of $X$ into $\PEX$) is quasiconvex in $(\PEX , d_{pel})$. In particular, 
the sets $L_\alpha$
are uniformly quasiconvex.

To prove hyperbolicity of $(\PEX , d_{pel})$, it suffices to prove that for all $K \geq 1$, there exists
$C = C(K)$ such that for all $a, b \in \PEX$,
$K$-quasigeodesic bigons in $(\PEX , d_{pel})$ joining $a, b$ are $C$-thin.
By the previous paragraph, we can assume without loss of generality, that the pair of quasigeodesics
in question lie in $Z_{a,b} \subset \PEX$. By strong relative hyperbolicity of 
$(\PEX,d_{pel})$ the two quasigeodesics enter and leave each element $L_\alpha \subset Z_{a,b}$
at nearby points (nearness dictated by the constant $K$). Since each $L_\alpha$ is hyperbolic (by definition
of partial electrocution) the two quasigeodesics track each other
inside each element $L_\alpha \subset Z_{a,b}$ they meet. Since they track each other off the horosphere-like sets
$L_\alpha$ this proves the hyperbolicity of $(\PEX , d_{pel})$ and the proof of Lemma \ref{pel}.
Lemma \ref{pel-track} now follows by applying Lemma \ref{farb1}. 
$\Box$

\subsection{Trees of Spaces}

\begin{defn} \label{tree}
(Bestvina-Feighn \cite{BF}) $P: X \rightarrow T$ is said to be a  tree of geodesic
 metric spaces satisfying
the q(uasi) i(sometrically) embedded condition if the geodesic
 metric space $(X,d)$
admits a map $P : X \rightarrow T$ onto a simplicial tree $T$, such
that there exist $ \epsilon$ and $K > 0$ satisfying the
following: \\
1) For all vertices $v\in{T}$,
$X_v = P^{-1}(v) \subset X$ with the induced path metric $d_{X_v}$ is
 a geodesic metric space $X_v$. Further, the
inclusions ${i_v}:{X_v}\rightarrow{X}$
are uniformly proper, i.e. for all $M > 0$, $v\in{T}$ and $x, y\in{X_v}$,
there exists $N > 0$ such that $d({i_v}(x),{i_v}(y)) \leq M$ implies
${d_{X_v}}(x,y) \leq N$.\\
2) Let $e = [0,1]$ be an edge of $T$ with initial and final vertices $v_1$ and
$v_2$ respectively.
Let $X_e$ be the pre-image under $P$ of the mid-point of  $e$.
 There exist continuous maps ${f_e}:{X_e}{\times}[0,1]\rightarrow{X}$, such that
$f_e{|}_{{X_e}{\times}(0,1)}$ is an isometry onto the pre-image of the
interior of $e$ equipped with the path metric. Further, $f_e$ is fiber-preserving,
i.e. projection to the second co-ordinate in ${X_e}{\times}[0,1]$ corresponds via $f_e$
to projection to the tree $P: X \rightarrow T$.\\
3) ${f_e}|_{{X_e}{\times}\{{0}\}}$ and
${f_e}|_{{X_e}{\times}\{{1}\}}$ are $(K,{\epsilon})$-quasi-isometric
embeddings into $X_{v_1}$ and $X_{v_2}$ respectively.
${f_e}|_{{X_e}{\times}\{{0}\}}$ and
${f_e}|_{{X_e}{\times}\{{1}\}}$ will occasionally be referred to as
$f_{e,v_1}$ and $f_{e,v_2}$ respectively. 
\end{defn}

A tree of spaces as in Definition \ref{tree} above is said to be a tree of
hyperbolic metric
spaces, if there exists $\delta > 0$ such that $X_v, X_e$ are all
$\delta$-hyperbolic for all vertices $v$ and edges $e$ of $T$.

\begin{defn} A tree $P: X \rightarrow T$ of geodesic
 metric spaces is said to be a  tree of relatively hyperbolic metric spaces
if in addition \\
4) each  vertex space $X_v$ is strongly hyperbolic relative to a collection of subsets $\HH_v$ and  each 
edge space $X_e$ is strongly hyperbolic relative to a collection of subsets $\HH_e$.  The individual sets
$H_{v,\alpha}\in \HH_{v}$ or $H_{e,\alpha}\in \HH_{e}$ will be called {\bf horosphere-like sets}.\\
5)  the maps $f_{e,v_i}$ above ($i = 1, 2$) are {\bf
  strictly type-preserving}, i.e. $f_{e,v_i}^{-1}(H_{v_i,\alpha})$, $i =
  1, 2$ (for
  any $H_{v_i,\alpha}\in \HH_{v_i}$)
 is
  either empty or some $H_{e,\beta}\in \HH_{e}$. Also, for all
$H_{e,\beta}\in \HH_{e}$, there exists $v$ and
$H_{v,\alpha}$, such that $f_{e,v} ( H_{e,\beta}) \subset H_{v,\alpha} $.\\
6) There exists $\delta > 0$ such that each $\EE (X_v, \HH_v )$ is $\delta$-hyperbolic.\\
7) The induced maps (see below) of the coned-off edge spaces into the
  coned-off vertex spaces $\hhat{f_{e,v_i}} : \EE ({X_e}, \HH_e ) \rightarrow
  \EE ({X_{v_i}}, \HH_{v_i})$ ($i = 1, 2$) are uniform quasi-isometries. This is called the
 {\bf qi-preserving electrocution condition}
 \end{defn}

Given the tree of spaces with vertex spaces $X_v$ and edge spaces $X_e$ there exists a naturally associated  tree whose vertex spaces are
$\EE (X_v, {\HH}_v)$ and edge spaces are
$\EE (X_e, {\HH}_e)$ obtained by simply coning off the respective horosphere like sets.
Condition (4) of the above definition ensures that we have natural inclusion maps of edge spaces
$\EE (X_e, {\HH}_e)$ into adjacent vertex spaces $\EE (X_v, {\HH}_v)$.

The resulting tree of coned-off spaces $P: \TC (X) \rightarrow T$
will be called the {\bf induced
  tree of coned-off spaces}. The resulting space will thus be denoted as
  $\TC (X)$ when thought of as a tree of spaces.

{\bf Remark:}  Strictly speaking, the induced tree
  exists
for any collection of vertex and edge spaces satisfying the {\em
  strictly type-preserving condition}. Hyperbolicity is not essential
  for the existence of the induced tree of spaces.

The {\bf cone locus} of $\TC (X)$,
(the induced tree of coned-off spaces), is
the graph (in fact a forest) whose vertex set $\VV$ consists of the
cone-points $c_v$ in the vertex set and whose edge-set $\EE$
consists of the
cone-points $c_e$ in the edge set. Thus the cone locus consists of
edges $c_e \times [0,1]$ with $c_e \times \{ 0 \}$ 
and $c_e \times \{ 1 \}$ identified with the appropriate $c_v$'s.
 The incidence relations are dictated by
 the incidence relations in $T$. To see that the cone locus
is a forest, note that a single edge space cannot have
more than one horosphere-like set mapping to
a common horosphere-like set in a vertex-set. Hence there are no induced loops in the
cone locus, i.e. it is a forest.

Note that connected components of the cone-locus can be naturally
identified with sub-trees of $T$. This is because \\
a) each edge and vertex of the cone locus is (respectively) an
edge and vertex of $T$. \\
b) There are no loops in the cone locus as otherwise two horosphere-like
sets in the same edge space would have to be attached to the same horosphere-like set in a vertex space contradicting
Condition (4) above. \\
 Each such connected component of the
cone-locus will be called a {\bf maximal cone-subtree}. The collection
of {\em maximal cone-subtrees} will be denoted by $\TT$ and elements
of $\TT$ will be denoted as $T_\alpha$. Further, each maximal
cone-subtree $T_\alpha$ naturally gives rise to a tree $T_\alpha$ of
horosphere-like subsets depending on which cone-points arise as
vertices and edges of $T_\alpha$. The metric space that $T_\alpha$
gives rise to will be denoted as $C_\alpha$ and will be referred to as
a {\bf maximal cone-subtree of horosphere-like spaces}. $g_\alpha : C_\alpha
\to T_\alpha$ will denote the induced tree of horosphere-like sets. $\GG$
will denote the collection of these maps.
The collection
of $C_\alpha$'s will be denoted as $\CC$. \\
{\bf Note:} Each $T_\alpha$ thus appears in two guises:\\
1) as a subset of $\TC (X)$ \\
2) as the underlying tree of $C_\alpha$\\

\smallskip

We shall have need for both these interpretations.

Since the natural inclusion $i_v\colon (X_v,{\HH}_v)\to (X,\CC)$ takes a horosphere-like set
$H_{v,\alpha}$ to a horosphere-like set $C_\alpha$ and the image of no two horosphere-like sets
in $X_v$ lie in the same horosphere-like set
$C_\alpha$,
$i_v$ will induce an embedding $\hat i_v\colon \EE (X_v, \HH_v) \to \TC (X)$. Let $\widehat X_v=\EE (X_v, \HH_v)$.
Using the fact that all edge-to-vertex space inclusions are strictly type-preserving, it follows that the induced maps $\hat i_v\colon \widehat X_v \to \TC (X)$
 are uniformly proper embeddings, that is,  for all $M>0$, 
 there exists $N>0$ such that for any $v\in T$ and $x,y\in \widehat X_v$,
 $d_{{\TC}(X)}(\hat i_v(x),\hat i_v(y))\leq M$ implies
 $d_{\widehat X_v}(x,y)\leq N$. This fact shall be useful later.

In \cite{mahan-reeves} conditions on trees of relatively hyperbolic metric spaces were found ensuring the hyperbolicity of $X$ relative to the collection $\CC$. For brevity
let us denote $\PP \EE (X, \CC, \GG , \TT ) = X_{pel}$

\smallskip

\noindent {\bf Remark:} \\
 If $X$ is a tree  of relatively hyperbolic spaces and $\CC$ is the collection of maximal cone subtrees of  horosphere like spaces $C_\alpha$, then the tree
 of coned-off spaces $P:
\TC(X) \to T$ can be thought of as obtained from $X$ by partially electrocuting each $C_\alpha$ to the cone subtree $T_\alpha$.

\smallskip

From lemma \ref{pel}, it follows that
\begin{lemma}
If $X$ is  hyperbolic relative to the collection $\CC$, then
$( X_{pel},d_{pel})$ is a hyperbolic metric space.
\end{lemma}

Thus we can treat the tree $P: \TC (X) \rightarrow T$
of coned-off spaces as a partially electrocuted space
$\PP\EE (X,\CC,\GG,\TT)=(X_{pel},d_{pel})$.
Note that
 $g_\alpha \colon C_\alpha\to T_\alpha$ collapses
 $C_\alpha$, the tree of horosphere-like spaces to the underlying tree $T_\alpha$.

\subsection{Preliminaries on Cannon-Thurston Maps}

For a hyperbolic metric space $X$, the Gromov bordification will be denoted by $\bbar X$.

\begin{defn}
Let $X$ and $Y$ be hyperbolic metric spaces and
$i : Y \rightarrow X$ be an embedding.
 A {\bf Cannon-Thurston map} $\bbar{i}$  from $\bbar{Y}$ to
 $\bbar{X}$ is a continuous extension of $i$ to the Gromov bordifications $\bbar{X}$ and $\bbar{Y}$. \end{defn}

 The following lemma, given in  \cite{mitra-ct}, gives a necessary and sufficient condition for the existence of Cannon-Thurston maps.

\begin{lemma}\cite{mitra-ct}\label{mit-cond}
 A Cannon-Thurston map $\bbar{i}$ from $\bbar Y $ to $\bbar{X}$ exists for the proper embedding $i\colon Y\to X$ if and only if
there exists a non-negative function $M(N)$ with $M(N)\rightarrow \infty$ as $N\rightarrow \infty$ such that the following holds: \\
Given $y_0 \in Y$, for all geodesic segments $\lambda$ in
$Y$ lying outside an $N$-ball around $y_0$ $\in Y$ any geodesic segment in $X$ joining the end points of $i(\lambda )$ lies outside the
$M(N)$-ball around $i(y_0)\in X$.
\end{lemma}

Note that due to stability of quasigeodesics, the above statement is also true if geodesics are replaced by quasigeodesics.

\medskip

Let $X$ and $Y$ be    hyperbolic relative to the collections $\HH_X$ and $\HH_Y$ respectively. Let $i\colon Y\to X$ be a strictly type-preserving proper embedding, i.e. for $H_Y\in \HH_Y$ there exists $H_X\in \HH_X$ such that $i(H_Y)\subset H_X$ and images of distinct horospheres-like sets in $Y$   lie in distinct horosphere-like sets in $X$. It follows easily from the fact that the inclusion of $H$ into $H^h$ is uniformly proper for all $H \in \HH_X$ or $\HH_Y$ that
the proper embedding $i\colon Y\to X$  induces a proper embedding  $i_h\colon \GG(Y,\HH_Y)\to \GG(X,\HH_X)$.

\begin{defn}
A Cannon-Thurston map is said to exist for the
 pair $X, Y$ of relatively hyperbolic metric spaces and a strictly type-preserving inclusion
$i: Y \rightarrow X$ if a Cannon-Thurston map  exists for the induced map $i_h\colon \GG(Y,\HH_Y)\to \GG(X,\HH_X)$ between the respective hyperbolic cones.
\end{defn}

We now give a criterion for the existence of Cannon-Thurston maps for  relatively hyperbolic spaces.
Let $Y, X$ be hyperbolic rel. $\YY , \XX$ respectively. Let $Y^h = \GG (Y, \YY ),
\widehat Y = \EE (Y, \YY )$ and $X^h = \GG (X, \XX ),
\widehat X = \EE (X, \XX )$.
Recall that $B^h_R(Z) \subset X^h$ denotes the $R$-neighborhood of $Z$ in $(X^h, d_h)$.
Also by the Definition of $\GXH$, recall that  distances in $(X^h, d_h)$ are proper functions
of distances in $(X, d)$. Similarly for $Y$ and $Y^h$.

\begin{lemma}\label{crit-relhyp}
A Cannon-Thurston map for $i\colon Y \to X$ exists if and only if
there exists a non-negative function $M(N)$ with
$M(N)\rightarrow \infty$ as $N\rightarrow \infty$ such that the
following holds: \\
Suppose $y_0\in Y$,  and $\hat \lambda$ in $\widehat Y$ is
an electric quasigeodesic segment  starting and ending  outside horospheres.
If  $\lambda^b = \hat \lambda \setminus \bigcup_{K \in \YY} K$
 lies outside an $B_N (y_0) \subset Y$,
then for any electric quasigeodesic $\hat \beta$ joining the
end points of $\hat i (\hat \lambda)$ in $\widehat X$,
$\beta ^b = \hat \beta \setminus \bigcup_{H \in \XX} H$ lies outside
 $B_{M(N)} (i(y_0)) \subset X$.
\end{lemma}

\noindent {\bf Proof:}  By Lemma \ref{mit-cond}, we note that a Cannon-Thurston map exists for the pair $(X^h, Y^h)$ if and only if there exists a non-negative function $M(N)$ with
$M(N)\rightarrow \infty$ as $N\rightarrow \infty$ such that the
following holds: \\
If a  geodesic $\lambda^h \subset Y^h$ lies
 outside $B_N^h (y_0) \subset Y^h$, then any geodesic $\beta^h \subset X^h$
joining the end-points of $i (\lambda^h)$
 lies outside $B_{M(N)}^h (i(y_0)) \subset X^h$.
 Now, from Lemma \ref{farb1}, and the fact that distances in $(X^h, d_h)$ (resp.
$(Y^h, d_h)$) are proper functions
of distances in $(X, d)$ (resp.
$(Y, d)$)
it follows that
$\lambda^h \subset Y^h$ (resp. $\beta^h \subset X^h$) lies outside
$B_N^h (y_0) \subset Y^h$ (resp. $B_{M(N)}^h (i(y_0)) \subset X^h$)
 if and only if
$\lambda^b $ (resp.
$\beta^b$) lies outside $B_{N_1(N)}(y_0) \subset Y$
(resp. $B_{M_1(N)}(i(y_0)) \subset X$), where $M(N), M_1(N), N_1(N) $ are
all proper functions of $N$ and each other. The result follows.
$\Box$

\medskip

Finally, we specialize Lemma \ref{crit-relhyp} to the case we are interested
in, viz. trees of relatively hyperbolic spaces.
Recall that $P: X\to T$ is a tree  of relatively hyperbolic spaces,
with vertex spaces
${X_v}$ and edge spaces ${X_e}$. $\HH_{v}$ and $\HH_{e}$
are the collections of horosphere-like sets in ${X_v}$ and  ${X_e}$ respectively.
 $P: \TC (X) \to T$ is the induced
 tree  of coned-off hyperbolic metric spaces, with vertex spaces
$\EE (X_v, \HH_v ) = \hhat{X_v}$ and edge spaces
$\EE (X_e, \HH_e ) =\hhat{X_e}$.
 $\TT$ is the collection of maximal cone subtrees $T_\alpha$ in $\TC (X)$ and $\CC$ is the
collection of maximal cone-subtrees of horosphere-like spaces $C_\alpha$ in $X$.
Note that $\EE{ (\TC (X), \TT )}$
 may be identified with $\EE (X, \CC )$. 
 Suppose that $a, b \in X \setminus \bigcup_{C_\alpha\in\CC} C_\alpha$. Let \\
$\bullet$ $\EE \EE (X) = \EE{ (\TC (X), \TT )} = \EE (X, \CC )$. (The double $\EE$ indicates
partial electrocution followed by electrocution.) \\
$\bullet$ $\GG (X, \CC ) = X^h$. \\
$\bullet$ ${\hat{\beta}}$ be an electric geodesic in $\EE \EE (X)$ \\
$\bullet$ $\beta$ be the relative geodesic in $X$ with the same underlying subset \\
$\bullet$ $\beta^h$ denote the  geodesic in $X^h$ joining its end-points $a,b \in X^h$ \\
$\bullet$ $\beta_{pel}$ be a (partially electrocuted) geodesic in $\TC (X)=X_{pel}$ joining $a, b$.

\begin{lemma}\label{crit-tree}
A Cannon-Thurston map for $i\colon X_v \to X$ exists if and only if
there exists a non-negative function $M(N)$ with
$M(N)\rightarrow \infty$ as $N\rightarrow \infty$ such that the
following holds: \\
Given $y_0\in X_v$,  and an electric quasigeodesic segment $\hat{\lambda}$ in $ \hhat{X_{v}}$ with end-points $a, b$  outside horospheres, if $\hat{\lambda} \setminus \bigcup_{H_{v,\alpha}\in \HH_v} H_{v,\alpha}$
lies outside $B_N(y_0) = \{ y \in  X_v: d_v(y,y_0) \leq N \}$, then for any
 partially electrocuted quasigeodesic $ \beta_{pel}$ joining $a, b$
 in $\TC ( X)$, $\beta_{pel}^b =  \beta_{pel} \setminus \bigcup_{C_\alpha \in \CC}
 C_\alpha$ lies outside $B_{M(N)} (i(y_0)) \subset X$, where $B_{M(N)} (i(y_0))$ denotes the $M(N)$-ball about $i(y_0)$ in the (non-electrocuted) metric $d$ on $X$.
\end{lemma}

\noindent {\bf Proof:} By Lemma \ref{farb1} and Lemma \ref{pel-track} we find that
there exists $D \geq 0$ such that
${\hat{\beta}}  \setminus \bigcup_{ C_\alpha\in \CC}C_\alpha$,
$\beta_{pel}  \setminus \bigcup_{ C_\alpha\in \CC} C_\alpha$,
${\beta}  \setminus \bigcup_{ C_\alpha\in \CC} C_\alpha$ and
${\beta^h}  \setminus \bigcup_{ C_\alpha\in \CC} C_{\alpha}^ h$
lie in a $D$ neighborhood of each other. The Lemma then
follows from Lemma \ref{crit-relhyp}. $\Box$

\medskip

Thus, to prove the existence of a Cannon-Thurston map, we have to show that
given a geodesic $\lambda^h$ lying outside large balls in $X_{v}^h$, we can
construct a partially electrocuted quasigeodesic $\beta_{pel}$ in $\TC (X)$ satisfying the condition that it lies outside large balls in $\TC ({X})$ away from the sets $C_\alpha$.

\section{Existence of Cannon-Thurston Maps}
\label{sec:ct}

Throughout this section, we will assume that trees of relatively hyperbolic spaces are as in Definition \ref{tree}.
We will also assume that  horosphere-like sets are uniformly  separated. Notation is carried over from the end of the previous section.

\noindent {\bf Sketch of Proof:} \\
As in \cite{mitra-trees}, the key step for proving the existence of a Cannon-Thurston map is to construct a hyperbolic ladder $B_{\hat\lambda}$ in $\TC ({X})$,  where $\hat{\lambda}$ is an electric
 geodesic segment in ${\widehat{X}}_{v_0}$
for some $v_0\in T$; and a large-scale Lipschitz retraction
 $\hat\Pi _{\hat\lambda}$ from  $\TC ({X})$ onto $B_{\hat\lambda}$. This proves the quasiconvexity of $B_{\hat\lambda}$.
Further, we shall show that if the underlying relative geodesic $\lambda$ of $\hat{\lambda}$ lies outside a large ball in $(X_{v_0},d_{X_{v_0}})$ modulo horospheres then $B_{\hat\lambda}$
lies outside a large ball in $X$ modulo
horospheres. Quasiconvexity of $B_{\hat\lambda}$  ensures that geodesics joining points
on $B_{\hat\lambda}$ lie close to it modulo horospheres.

 We consider here electric geodesics $\hat{\mu}$
in the coned-off vertex and edge-spaces $\hhat{X_v}$ and $\hhat{X_e}$.
In \cite{mitra-trees}, we have assumed that each $X_v, X_e$ are $\delta$-hyperbolic metric spaces and took $\lambda=\hat\lambda$, hence we needed to find points in some $C$-neighborhood of $\lambda$ to construct $B_\lambda$. Since there is only the usual (Gromov)-hyperbolic metric in \cite{mitra-trees}, this creates no confusion. But, in the present situation,
we have two metrics $d_{X_v}$ and $d_{\widehat X_v}$ on $X_v$. As electrically close (in the $d_{\widehat X_v}$ metric)
does not  imply close (in the $d_{X_v}$ metric),
we cannot take a $C$-neighborhood in the $d_{\widehat X_v}$ metric.
 Instead  we will first construct  an electroambient representative
$\lambda $ of $\hat \lambda$ in the space $X_{v}^h$ and
take a  neighborhood of $\lambda $ in  $X_{v}^h$.
Also, by noting that $X_{v}^h$ with {\em horoballs} coned off is the same as $\hhat{X_v}$, we will be able to carry out the construction in
\cite{mitra-trees} mutatis mutandis.

Finally we construct vertical quasigeodesic
rays in $B_{\hat\lambda}$ to show that if
$\hat\lambda \setminus\bigcup_{H_{v\alpha}\in\HH_v}H_{v \alpha}$
lies outside a large ball in $X_v$, then
$(B_{\hat\lambda}\setminus \bigcup_{C_\alpha\in\CC} C_{\alpha})$ lies outside a large ball in $X$. The existence of a Cannon-Thurston map follows from Lemma \ref{crit-tree}.

\subsection{Construction of Hyperbolic Ladder}
\label{subsec:cons}

Given a geodesic segment $\hat \lambda \subset \widehat X_{v_0}$ with end points lying outside horospheres-like sets,  we now construct a quasiconvex set $B_{\hat \lambda} \subset \widehat X$ containing $\hat \lambda$.

The quasi-isometric embeddings $f_{e,v}\colon X_e\to X_v$ induce quasi-isometric embeddings $f_{e,v}^h\colon X_{e}^h\to X_{v}^h$. Thus for every
 vertex $v$ and edge $e$ in $T$, $f_{e,v}^h(X_{e}^h)$ will be
$C_2$ quasiconvex for some $C_2>0$ as $X_v^h$ is a uniformly
hyperbolic metric space for all $v$. Let $\delta$ be the
hyperbolicity constant of  $X_{v}^h, X_{e}^h$. Let $C_1$ be as in Lemma \ref{gen lem}.
Let $C=C_1+C_2$

For $Z\subset X_{v}^h$, let $N_C(Z)$ denote the $C$-neighborhood of $Z$ in $X_{v}^h$,  where $C$ is as above.

\smallskip

\noindent {\textbf{Hyperbolic Ladder $B_{\hat\lambda}$}}

\smallskip

 Recall that $P: \TC (X) \rightarrow T$ is the usual projection to the base tree.

For convenience of exposition, $T$ shall be assumed to be rooted, i.e.
equipped with a base vertex $v_0$. Let $v \neq v_0$ be a vertex of $T$.
Let $v_{-}$ be the penultimate vertex on the geodesic edge path
from $v_0$ to $v$. Let $e$ denote the directed edge from ${v_{-}}$
to $v$.

Define  $\phi _{v,e}\colon f_{e,v_{-}}(X_e)\to  f_{e,v}(X_e )$ as follows: \\
If $p\in  f_{e,v_{-}}(X_e) \subset X_{v_{-}}$, choose $x\in X_e$ such that
$p=f_{e,v_{-}}(x)$ and define   $\phi _{v,e}(p)= f_{e,v}(x )$.

Note that in the above definition, $x$ is chosen from a set of bounded diameter.

\smallskip

Since  $f_{e,v_{-}}$ and $f_{e,v}$ are quasi-isometric embeddings into their respective vertex spaces $\phi_{v,e}$'s are uniform quasi-isometries for all vertices.
We shall denote $\EE (X_v , \HH_v) = \EE ( \GG (X_v , \HH_v), \HH_v^h ) = \hhat{X_v}$ and
$\EE (X_e , \HH_e) = \EE ( \GG (X_e , \HH_e), \HH_e^h ) = \hhat{X_e}$.

\smallskip

\noindent \underline{Step 1}\\ Without loss of generality, assume that the base tree $T$ is rooted with base vertex $v_0$. Let $v$ be any vertex of $T$. 
Let $\hat \mu \subset \hhat{X_v}$ be a geodesic segment in
$(\hhat{X_v}, d_{\hhat{X_v}} )$ with starting and ending points lying outside horoballs and $\mu$ be the corresponding electroambient quasi-geodesic in $X_{v}^h$ (cf Lemma \ref{ea}).   Then $P(\hat \mu)=v$.
For the collection of edges $e^\prime$ in $T$ incident on $v$ , but not lying on the geodesic (in T) from $v_0$ to $v$,
consider the subcollection of edges $\{ e\}$ for which $N_C^h(\mu) \cap f_{e,v}(X_e) \neq \emptyset$
and for each such $e$,
 choose $p_e, q_e\in N_C^h(\mu) \cap f_{e,v}(X_e)$ such that
$d_{X_{v}^h}(p_e, q_e)$ is maximal. Let
$\{v{_i}\}$  be the terminal vertices of edges $e_i$ for which
$d_{\hhat X_{v}}(p_{e_i}, q_{e_i})> D$ ,
where $D$ is as in Lemma \ref{gen lem} above. Let $\hat\mu_{v,e_i}$ be a geodesic in $\hhat X_v$ joining $\phi_{v_i,e_i}(p_{e_i})$ and
$\phi_{v_i,e_i}(q_{e_i})$.
Define $$B^1(\hat \mu)= \hat{i_v}(\hat \mu)\cup \bigcup_{i}
\hat \Phi _{v,e_i}(\hat \mu_{v,e_i}) $$
where $\hat \Phi _{v_i,e_i}(\hat \mu_{v,e_i})$ is an electric geodesic in $\widehat X_{v_i}$
joining $\phi_{v_i,e_i}(p_{e_i})$ and $\phi_{v_i,e_i}(q_{e_i})$.

\smallskip

\noindent \underline{Step 2} \\
Step 1 above constructs $B^1(\hat \lambda)$ in particular.   We proceed inductively.  \\
 Suppose that $B^m(\hat \lambda)$ has been constructed such that the vertices
in
 of $P(B^m(\hat \lambda))\subset T$ are the vertices of a subtree.  Let $\{w_i\}_{i} = P(B^m(\hat \lambda))\setminus P(B^{m-1}(\hat \lambda))$. \\
Assume further that $P^{-1}(w_k)\cap B^m(\hat \lambda)$ is a path of the form $i_{w_k}(\hat \lambda_{w_k})$ , where $\hat \lambda_{w_k}$ is a geodesic in $(\widehat X_{w_k}, d_{\widehat X_{w_k}})$ .  \\
Define
$$B^{m+1}(\hat \lambda)= B^m(\hat \lambda)\cup \bigcup _{k} (B^1(\hat \lambda_{w_k}))$$
where  $B^1(\hat \lambda_{w_k})$ is defined in step 1  above.  \\
Define
$$B_{\hat \lambda} =\cup _{m\geq 1}B^m(\hat \lambda)$$
Observe that the vertices comprising $P(B_{\hat \lambda})$ in $T$ are
the vertices of a subtree,  say, $T_1$.

\smallskip

Roughly speaking, what we have done is that at each stage we take a geodesic $\hat \lambda_v$
look at all edge spaces which hit $\widehat X_v$ near $\hat \lambda_v$, `break'
$\hat \lambda_v$ into maximal subpieces coarsely contained
in the images of these edge spaces and then `flow' them (via the $[0,1]$ direction in $\widehat X_e \times
[0,1]$) into adjacent vertex spaces. The maximal subpieces are the $\hat \mu$'s.

\subsection{Retraction Map}

In order to prove $B_{\hat\lambda}$ is quasiconvex in $\TC(X)$, we will construct a \textit{retraction map} $\hat\Pi_{\hat\lambda}$ from $\TC(X)$ to $B_{\hat\lambda}$ which is coarsely Lipschitz. For convenience of exposition, we shall define $\hat\Pi_{\hat\lambda}$ only on the union of vertex spaces of $\TC(X)$.\\

For a tree $T$, let $\VV(T)$ denote the vertex set of $T$.\\
\begin{defn} \label{retdef}
Let $ \hat \pi_{\hat \lambda_v}\colon \widehat X_v\to \hat \lambda_v$ be an electric projection from $\widehat X_v$ onto $\hat \lambda _v$ (See Definition \ref{proj}).

Define $\widehat \Pi_{\hat \lambda}$  on $\bigcup _{v\in\VV(T_1)}\widehat X_v$ by
\begin{center}
$\widehat \Pi_{\hat \lambda} (x)= \hat i_v(\hat \pi_{\hat \lambda_v}(x))$ for $x\in \widehat X_v$.
\end{center}
\medskip
If $x\in $ $P^{-1}((\VV(T)\setminus \VV(T_1)))$,  choose $x_1\in$ $ P^{-1}(\VV(T_1))$ such that $d(x, x_1)=d(x, P^{-1}(\VV(T_1)))$ and define $\widehat \Pi ' _{\hat \lambda} (x)=x_1$. Next define $\widehat \Pi_{\hat \lambda}(x)=\widehat \Pi_{\hat \lambda}(\widehat \Pi ' _ {\hat \lambda}(x))$.
\end{defn}
\begin{theorem} \label{rel-ret thm}
There exists $C_0 \geq 0$ such that

\begin{center}
$d_{\TC(X)}(\widehat \Pi_{\hat \lambda}(x), \widehat \Pi_{\hat \lambda}(y))\leq C_0 d_{\TC(X)}(x, y)+C_0$ for $x, y\in \TC(X)$.
\end{center}
In particular, if $\TC(X)$ is hyperbolic, then $B_{\hat \lambda}$ is uniformly (independent of $\hat\lambda$) quasiconvex.
\end{theorem}

Though we have made some small changes from \cite{mitra-trees} in the construction of $B_{\hat \lambda}$,  the proof of the above theorem will be quite similar to the proof of Theorem 3.8 in \cite{mitra-trees}. We include the proof here for completeness, as there are slight differences at various stages.

Also, we state a couple of
 lemmas (proven in \cite{mitra-trees}) which will be required to prove Theorem \ref{rel-ret thm}.
In the construction of the
hyperbolic ladder $B_{\hat\lambda}$, recall that we had `broken' the
 quasigeodesics $\lambda_v$  into `maximal subpieces' lying close to 
 $\lambda_v$
and having end points lying in the same edge space.

Lemma \ref{near lem} below
 says that points in the corresponding edge space
 which are at bounded distance from $\lambda_v$'s are also at bounded distance from `maximal subpieces'.

Recall from the construction of the hyperbolic ladder $B_\lambda$
that $f^h_{e,v}(X_e^h)$ are uniformly $C_2$-quasiconvex and $C=C_1+C_2$,
where $C_1$ is as in Lemma \ref{gen lem}.
\begin{lemma}\label{near lem} (See Lemma 3.6 of \cite{mitra-trees})
Let $\hat \mu_1 \subset \widehat X_v$ be an  electric geodesic with end points $a$ and $b$ lying outside horosphere.   Let $\mu _1$ be the corresponding electroambient quasigeodesic in $X_v^h$.   Let $p, q\in N_C(\mu_1)\cap f_{e,v}^h(X_e^h)$ be such that $d_{X_v^h}(p, q)$ is maximal.  Let $\hat \mu_2$ be a geodesic in $\widehat X_{v}$ joining $p$ and $q$ and $\mu_2$ be its electro-ambient representative.  If $r\in N_C(\mu_1)\cap f_{e,v}^h(X_e^h)$ then $d_{X_v^h}(r, \mu_2)\leq P_6$ for some constant $P_6$ depending only on $C, D, \delta$.
\end{lemma}

Lemma \ref{pi lem}
below says that images of points in edge spaces  under
nearest point projections to $\lambda_v$'s and `maximal subpieces'
  are at a bounded  distance from each other.
\begin{lemma}\label{pi lem}(See Lemma 3.7 of \cite{mitra-trees})
Let $\hat \mu_1 $ and $\hat \mu_2$ be as in Lemma \ref{near lem}. If $z\in f_{e,v}(X_e)$,  then $d_{X_v^h}(\pi_{\mu_1}(z),\pi_{\mu_2}(z))\leq P_7$ for some constant $P_7$ depending only on $\delta, C, D$.
In particular, if $z\in \hhat f_{e,v}(\widehat X_e)$, then $d_{\widehat X_v}(\hat \pi_{\hat \mu_1}(z), \hat \pi_{\hat \mu_2}(z))\leq P_7$.	
\end{lemma}

Suppose $x,y\in{\widehat X_v} \subset
\TC(X)$ and $d_{\TC(X)}(x,y)\leq 1$. Since $\hat i_v$'s are uniformly proper embedding from $\hhat X_v$ to $\TC(X)$, there exists a constant $M>0$ such that $d_{\widehat X_v}(x,y)\leq M$.\\
Let ${\mathbb P} =max\{P_i,M:1\leq i\leq 7\}$, where $P_i$'s are the
constants
in Lemmas \ref{proj lem}, \ref{q.i. lem}, \ref{bdd disc}, \ref{ret lem},  \ref{e.qi}, \ref{near lem} and \ref{pi lem}.

\medskip

{\bf Proof of theorem \ref{rel-ret thm}}

It suffices to prove that if $d_{\TC(X)}(x, y)\leq 1$ then $d_{\TC(X)}(\widehat \Pi_{\hat \lambda}(x), \widehat \Pi_{\hat \lambda}(y))\leq C_0$.

Let $d_{\TC(X)}(x, y)\leq 1$.

\underline{\textbf{Case 1}}: Let $x, y\in P^{-1}(v)$ for some $v\in T_1$. Using Lemma \ref{ret lem}, there exists a constant $K_0({\Bbb P})>0$ such that
\begin{center}
$d_{\TC(X)}(\widehat \Pi _{\hat \lambda_v}(x),\widehat \Pi _{\hat \lambda_v}(y))
\leq d_{\widehat X_{v}}(\hat \pi _{\hat \lambda_v}(x),\hat \pi _{\hat \lambda_v}(y))\leq K_0$.                                       \end{center}

\medskip

\underline{\textbf{Case 2}}: Let $x\in P^{-1}(w)$ and $y\in P^{-1}(v)$  for some $v, w\in T_1$ such that $v\neq w$. Now $v$ and $w$ are adjacent in $T_1$ since $d_{\TC(X)}(x, y)\leq 1$.  Without loss of generality we can assume that $w=v_-$. Let $e$ be the edge between $v$ and $w$.\\
Recall that
\begin{center}
 $B_{\hat \lambda}\cap P^{-1}(v)=\hat \lambda_v$, \\$B_{\hat \lambda}\cap P^{-1}(w)=\hat \lambda_w$, \\
$\hat \lambda_v =\hat \Phi_{v,e}(\hat \mu_{w,e})$,
\end{center}
\noindent and end points of $\mu_{w,e}$ lie in a $C$-neighborhood of $\lambda_w$.

\noindent \underline{Step 1: }From lemma \ref{pi lem},
\begin{center}
$d_{\TC(X)}(\widehat \Pi_{\hat \lambda _w}(x), \widehat \Pi_{\hat \mu _w}(x))\leq d_{\widehat X_w}(\hat \pi_{\hat \lambda _w}(x), \hat \pi_{\hat \mu_w}(x))\leq {\Bbb P}$.                                                                                             \end{center}
\underline{Step 2:} $\hhat f_{e,v}(\hhat X_e)$ and $\hhat f_{e,w}(\hhat X_e)$ are uniformly quasiconvex in $\hhat X_v$ and $\hhat X_w$ respectively. Therefore, by stability of quasigeodesics and using Lemma \ref{e.qi}, there exists a constant $R>0$ such that
$d_{\widehat X_v}(\hat\phi_v(\hat \pi_{\hat \mu_w} (x)), \hat \pi_{\hat \lambda_v}(\hat\phi_v(x)))\leq R$.
\begin{eqnarray*}
d_{\TC(X)}(\widehat \Pi_{\hat \mu_w} (x), \widehat \Pi_{\hat \lambda_v}(\hat\phi_{v,e}(x)))&=&
d_{\TC(X)}(\hat \pi_{\hat \mu_w} (x), \hat \pi_{\hat \lambda_v}(\hat\phi_{v,e}(x)))\\
&\leq&d_{\TC(X)}(\hat \pi_{\hat \mu_w} (x), \hat\phi_{v,e}(\hat \pi_{\hat \mu_w} (x))+\\&&d_{\TC(X)}(\hat\phi_{v,e}(\hat \pi_{\hat \mu_w} (x)), \hat \pi_{\hat \lambda_v}(\hat\phi_{v,e}(x)))\\
&\leq& 1+R.
\end{eqnarray*}
\underline{Step 3:} $d_{\TC(X)}(x, y)=1=d_{\TC(X)}(x,\hat\phi_{v,e}(x))$. Then $d_{\widehat X_v}(\hat\phi_{v,e}(x), y)\leq 2M$.
Thus using lemma \ref{ret lem}, we have \\
$d_{\TC(X)}(\widehat \Pi_{\hat \lambda_v}(\hat\phi_{v,e}(x)), \widehat \Pi_{\hat \lambda_v}(y))\leq d_{\widehat X_v}(\hat \pi_{\hat \lambda_v}(\hat\phi_{v,e}(x)), \hat \pi_{\hat \lambda_v}(y))\leq 2{\Bbb P}M + {\Bbb P}$. \\
Thus from above three steps, there exists a constant $K_1({\Bbb P})>0$ such that
$d_{\TC(X)}(\hat\Pi_{\hat \lambda}(x), \hat\Pi_{\hat \lambda}(y))\leq K_1$.
\medskip

\underline{\textbf{Case 3}}:
$P([x, y])$ is not contained in $T_1$. Then $P(x)$ and $P(y)$ belong to the closure of the same component of $T\setminus T_1$. Then $P(\widehat \Pi '_{\hat \lambda}(x))=P(\widehat \Pi '_{\hat \lambda}(y))=v$ for some $v\in \VV(T_1)$ by the second part of Definition \ref{retdef}. 
Let $x_1=\widehat \Pi '_{\hat \lambda}(x)$ and $y_1=\widehat \Pi '_{\hat \lambda}(y)$.  Now $x_1,y_1 \in f_{e,v}^h(X_e^h)$ for some edge $e$ with initial vertex $v$. (This
is because for such
$x_1, y_1$,  Definition 
\ref{retdef} shows that $\Pi '_{\hat \lambda}$ is given as a two-step process: first choose $x
\in X_v$ for the vertex $v\in \VV(T_1)$ closest to $P(x_1)$ and
then project $x$. Since $x, y$ project onto the same component of $T\setminus T_1$ under $P$, the shortest paths from
$x_1, y_1$ to $P^{-1}(T_1)$ must enter $P^{-1}(T_1)$ through
the same edge space.)\\
If $d_{X_v^h}(\pi_{\lambda_v}(x_1), \pi_{\lambda_v}(y_1))<D$ then
\begin{center}
$d_{\TC(X)}(\widehat \Pi_{\hat \lambda}(x), \widehat \Pi_{\hat \lambda}(y))\leq d_{\widehat X_{v}}(\hat \pi _{\hat \lambda_v}(x),\hat \pi _{\hat \lambda_v}(y))\leq D$.                                                                                 \end{center}
Let us assume $d_{X_v^h}(\pi_{\lambda_v}(x_1), \pi_{\lambda_v}(y_1))>D$.
Then by Lemma \ref{gen lem}
$[x_1,\pi_{\lambda_v}(x_1)]\cup [\pi_{\lambda_v}(x_1), \pi_{\lambda_v}(y_1)]\cup [ \pi_{\lambda_v}(y_1),y_1]$ is a quasi-geodesic lying in a $C_1$-neighborhood of $[x_1, y_1]$. \\
Since $f_{e,v}^h(X_e^h)$ is $C_2$-quasiconvex in $X_v^h$, there exists $x_2,y_2\in f_{e,v}^h(X_e^h)$ such that $d_{X_v^h}(\pi_{\lambda_v}(x_1),x_2)\leq C_1+C_2=C$ and $d_{X_v^h}(\pi_{\lambda_v}(y_1),y_2)\leq C_1+C_2=C$. Thus $x_2,y_2\in N_C(\lambda_v)\cap f_{e,v}^h(X_e^h)$. There exists $x_2',y_2'\in N_C(\lambda_v)\cap f_{e,v}(X_e)$ such that $d_{\widehat X_v}(x_2,x_2')\leq 1$ and $d_{\widehat X_v}(y_2,y_2')\leq 1$.\\
Let $D_1>P_4D+P_4$. If $D_1<d_{\widehat X_v}(\hat \pi_{\hat \lambda_v}(x_2'), \hat \pi_{\hat \lambda_v}(y_2'))$ then by Lemma \ref{ret lem} $d_{\widehat X_v}(x_2',y_2')>D$. This implies the edge $P(e)$ of $T$ would be in $T_1$, (because we would be able to to continue the construction of the ladder $B_\lambda$ beyond the vertex $v$ through the edge $e$) 
 which is a contradiction.
Therefore $d_{\widehat X_v}(\hat \pi_{\hat \lambda_v}(x_2'), \hat \pi_{\hat \lambda_v}(y_2'))\leq D_1$.
\begin{eqnarray*}
d_{\TC(X)}(\widehat \Pi_{\hat \lambda}(x), \widehat \Pi_{\hat \lambda}(y))
&\leq& d_{\widehat X_v}(\hat \pi_{\hat \lambda_v}(x_1), \hat \pi_{\hat \lambda_v}(y_1))\\
&\leq& d_{X_v^h}(\pi_{\lambda_v}(x_1),x_2)+d_{X_v^h}(\pi_{\lambda_v}(y_1),y_2)\\&&+d_{\widehat X_v}(x_2',y_2')+2\\
&\leq& 2C+d_{\widehat X_v}(x_2',\hat \pi_{\hat \lambda_v}(x_2'))+d_{\widehat X_v}(y_2',\hat \pi_{\hat \lambda_v}(y_2'))\\&&+d_{\widehat X_v}(\hat \pi_{\hat \lambda_v}(x_2'), \hat \pi_{\hat \lambda_v}(y_2'))+2\\
&\leq& 4C+D_1+2=K_2(say).
\end{eqnarray*}

Taking $C_0 = $max$\{K_0, K_1, K_2,D\}$,  we have the required result. $\Box$

\subsection{Vertical Quasigeodesic Rays}
\label{subsec:qiray}
Let  $\hat\lambda _{v_0}$  be an electric geodesic segment from $a$ to $b$ in $\widehat X_{v_0}$ with $a$ and $b$  lying outside horosphere-like sets and $\lambda _{v_0}$  denotes its  electroambient quasigeodesic representative in $X_{v_0}^h$.

Note that
$B_{\hat \lambda _{v_0}}=\bigcup_{v\in\VV(T_1)}\hat i_v(\hat \lambda _{v})$
is the quasiconvex set constructed above, where $T_1$ is the subtree of $T$ defined above.

\smallskip

Let \\
$\bullet$ $\lambda ^c _{v}$ be the union of geodesic subsegments of the electroambient quasigeodesic $\lambda _{v}$ lying inside  the horoball-like sets penetrated by $\lambda _{v}$.  \\
$\bullet$
 $\lambda ^b _{v}=\lambda _{v}\setminus \lambda ^c _{v}$.
(Note that $\lambda _{v}^b\subset \hat \lambda _{v}$).  \\
$\bullet$
 $B^b_{\lambda _{v_0}}= \bigcup _{v\in T_1 }i_v(\lambda ^b _v)$.
(Then $B^b_{\lambda_{v_0}}\subset B_{\hat \lambda _{v_0}}$). \\

\smallskip

If $x\in B^b_{\lambda_{v_0}}$, then there exists $v\in T_1$ such that
$x\in\lambda^b_{v}$.
Let \\
$S = [v_n,v_{n-1}]\cup...\cup[v_1,v_0]$ be the geodesic edge path in $T_1$ joining $v$ and $v_0$.

We will construct a map
$r_x\colon S\to B^b_{\lambda_{v_0}}$ satisfying \\
$\bullet$ $d_{S}(w, w')\leq d(r_x(w), r_x(w')\leq C d_{S}(w, w')$ for all $w,w' \in S$, where $d$ denotes the metric on $X$.\\
$\bullet$ $r_x(v_i) \in X_{v_i}$. \\

\smallskip

$r_x$ will be called a \textit{ vertical quasigeodesic ray}.

\smallskip

\noindent $\bullet$ For $v_n\in S$, define $r_x(v_n)=x$

\smallskip

Let $v=v_n,w=v_{n-1}$, $e_i=[v_i,v_{i-1}]$ and $\psi_{e_i,v_i}=\phi_{v_{i}}^{-1}\colon f_{e_i,v_i}(X_{e_i})\to f_{e_{i},v_{i-1}}(X_{e_i})$ for all $i=1,...,n$.

Then $\psi_{e_i,v_i}$ is a quasi-isometry.

Since $x$ lies outside horosphere-like sets and $\psi_{e,v}$ preserves horosphere-like sets
(by the strictly type-preserving condition),
 $\psi_{e,v}(x)$ will lie outside horosphere-like sets.

Let $[a,b]$ be the maximal connected component of $\lambda_v ^b$ on which $x$ lies. Then there exists two horosphere-like sets $H_1$ and $H_2$ such that $a\in H_1$ (or is a initial point of $\lambda_v$) and $b\in H_2$ (or is a terminal point of $\lambda_v$). Since $\psi_v$ preserves horosphere-like sets, $\psi_v([a,b])\setminus\{\psi(a),\psi(b)\}$ will lie outside horosphere-like sets.

As $ \psi_v$ is a quasi-isometry, $ \psi_v([a,b])$ will be a quasi-geodesic in $X_{w}$. Let $ \Psi_v^h([a,b])$ be the hyperbolic geodesic in $X_{w}^h$ joining $ \psi_v(a)$ and $ \psi_v(b)$. Then $ \psi_v([a,b])$ will lie in a bounded neighborhood of $ \Psi_v^h([a,b])$. Therefore there exists $C_1>0$ such that $d( \psi_v(x_1), \Psi_v^h([a,b]))\leq C_1$.
By Lemma \ref{farb1}, there exists an upper bound on how much $ \Psi_v^h([a,b])$ can penetrate horoball-like sets, that is, for all $z\in  \Psi_v^h([a,b])$ there exists $z'\in  \Psi_v^h([a,b])$ lying outside horoball-like sets such that $d(z,z')\leq B$. Hence there exists $y_1\in  \Psi_v^h([a,b])$ such that $d(\psi_v(x_1),y_1)\leq B+C_1$ and $y_1$ lies outside horosphere-like sets.

Again, $ \Psi_v^h([a,b])$ lies at a uniformly bounded distance $\leq C_2$ from $\mu_v$ (the electroambient quasigeodesic representative of $\hat \mu_v$  in the construction of $B_{\hat\lambda_{v_0}}$)). Let $c,d\in \mu_v$ such that $d(a,c)\leq C_2$ and $d(b,d)\leq C_2$. Then $ \Psi_v^h([a,b])$ and the quasigeodesic segment $[c,d]\subset \mu_v$ have similar intersection patterns (Lemma \ref{farb1} ) with horoball-like sets.
Therefore $[c,d]$ can penetrate only a bounded distance $\leq B$ into any horoball-like set. Hence there exists $y_2\in \mu_v$ and $y_2$ lies outside horosphere-like sets such that $d(y_1,y_2)\leq C_2+B$.

Since end points of $\mu_v$ lie at a bounded neighborhood of $\lambda_w$, there exists $C_3>0$ such that $\mu_v$ will lie at a $C_3$ neighborhood of $\lambda_w$. Therefore there exists $y_3\in \lambda_w$ such that $d(y_2,y_3)\leq C_3$. Now $y_3$ may lie in a horoball-like set.
Since $\mu_v$ and $ \pi_{\lambda_w}(\mu_v)$ lies in a bounded neighborhood of each other, by Lemma \ref{farb1} they have similar intersection patterns with horoball-like sets. Therefore there exists $B>0$ and $y\in  \lambda_w$ such that $y$ lies outside horosphere-like sets and $d(y_3,y)\leq B$.

Hence $d(x,y)\leq 1+C_1+C_2+C_3+3B = C\mbox{(say)}$.

\smallskip

\noindent $\bullet$ Recall that $w=v_{n-1}$. Define $r_x(v_{n-1})=y$.  \\

\smallskip

Thus we have $1\leq d(r_x(v_n), r_x(v_{n-1}))\leq C$.

Using the above argument  repeatedly (inductively replacing $x$ with $r_x(v_i)$ in each step) we get the following. Since $r_x(v) \in X_v$, we have $d_S(v,w)\leq d(r_x(v), r_x(w))$.

\begin{lemma}\label{ray lem}
There exists $C\geq 0$ such that for all $x\in \lambda _v ^b \subset B_{\lambda _{v_0}} ^b $ there exists a $C$-vertical quasigeodesic ray $r_x\colon S \to B_{\lambda _{v_0}} ^b$ such that $r_x(v)=x$
and $d_S(v,w)\leq d(r_x(v), r_x(w)\leq C d_S(v, w)$, where $S$ is the geodesic edge  path in $T_1$ joining $v$ and $v_0$ and $w\in S$.
\end{lemma}

The following is the concluding Lemma of this subsection.

\begin{lemma}\label{final lem}
There exists $C>0$ such that the following holds.
Fix a  reference point
$p$ lying outside the horospheres in $X_{v_0}$. $B_n(p)$ denotes the $n$-ball
around $p$ in $(X_{v_0}, d_{X_{v_0}})$.
Let $x\in \lambda ^b _v \subset B^b _{\lambda _{v_0}}\subset B_{\hat \lambda _{v_0} }$ such
that
 $\lambda ^b _{v_0}$ lies outside $B_n(p)$ (and
hence
 entry and exit points of $\hat \lambda$ to a horosphere lie
outside $B_n(p)$).  Then $x$ lies outside an $n/(C+1)$ ball
about $p$ in $X$.
\end{lemma}

\noindent {\bf Proof:} Since
$\lambda ^b _{v_0}$ is a part of $\hat\lambda _{v_0}$,  $r_x(v_0)$ lies outside $B_n(p)$.

Let $m$ be the first non-negative integer such that
$v\in P(B^m(\hat \lambda _{v_0}))\setminus P(B^{m-1}(\hat \lambda _{v_0}))$.
Then $d_{T_1}(v_0, v) = m$,
and $d(x, p)\geq m$ (since $r_x(v)=x\in \lambda ^b _v$).\\
From Lemma \ref{ray lem}, $m\leq d(r_x(v), r_x(v_0))\leq C m$.  \\
Since $r_x(v_0)$ lies outside $B_n(p)$, $d(r_x(v_0), p)\geq n$. \\
$n\leq d(r_x(v_0), p)\leq d(r_x(v_0), r_x(v)) + d(r_x(v), p) \leq mC  + d(r_x(v), p)$.  \\
Therefore,  $d(r_x(v), p)\geq n-mC$ and $d(r_x(v), p)\geq m$.  \\
Hence $d(x, p)=d(r_x(v), p)\geq \frac {n}{C+1}$.
 $\Box$

\subsection{Proof of Main theorem}
\label{proof}
First we recall the following notation: \\
$\bullet$ $\hat\lambda _{v_0}$ = electric geodesic in $\hhat X_{v_0}$ joining $a, b$ with $\lambda ^b _{v_0}$ lying outside $B_n(p)$ for a fixed reference point $p\in X_{v_0}$ lying outside horosphere-like sets.  \\
$\bullet$ $\lambda _{v_0}$ = electroambient quasi-geodesic in $X_{v_0}^h$ constructed from $\hat \lambda _{v_0}$.  \\
$\bullet$ $\beta_{pel}$ = quasi-geodesic in the partially electrocuted space $ \TC ( X)=(X_{pel},d_{pel}).$
 joining $a,b$.\\
$\bullet$ $\beta $ = electroambient quasi-geodesic in $\GXH$ corresponding to $\beta_{pel}$.  \\
$\bullet$ $\beta'_{pel} = \hat \Pi_{\hat\lambda _{v_0}}(\beta_{pel})$,  where $\hat\Pi _ {\hat \lambda _{v_0}}$ is the retraction map from the vertex space of $\TC(X)$ to the quasi-convex set $B_{\hat \lambda _{v_0}}$. \\

By Lemma \ref{pel}, $ \TC ( X)$ is hyperbolic.
$\beta _{pel}'$ is a dotted  quasi-geodesic (i.e. a quasi-isometric
embedding of the integer points in an interval)
in the partially electrocuted space $ \TC ( X)$ lying on $B_{\hat \lambda _{v_0}}$.\\So $\beta _ {pel}'$ lies in a $K$-neighborhood of $\beta _{pel}$ in $ \TC ( X)$. But $\beta _{pel}'$ might backtrack. $\beta _{pel}'$ can be modified to form a  quasigeodesic in $ \TC ( X)$ of the same type (i.e. lying in a $K$-neighborhood of $\beta _{pel}$ in $ \TC ( X)$) without backtracking with end points remaining the same. Let $u$ and $v$ be the first and last points on a horosphere-like set $C_\alpha$ penetrated by $\beta _{pel}'$. Replace $\beta _{pel}'$ on $[u,v]_{  \TC ( X)} (\subset \beta _{pel}')$  by a partially electrocuted geodesic from $u$ to $v$ in $ \TC ( X)$. Since $\beta _{pel}'$ has finite length in $ \TC ( X)$, it jumps from one horosphere to another a finite number of times, so after a finite number of modification we will have produced a path that does not backtrack. In this process, we have if anything, reduced the length of $\beta _{pel}'$, so the new path, $\gamma _{pel}$ (say), is again a partially electrocuted quasigeodesic of the same type as $\beta _{pel}'$.

Now $a,b$ are end points of $\hat \lambda _{v_0}$, therefore $a,b\in B_{\hat \lambda _{v_0}}$ and end points of
$\gamma _{pel}$ are $a,b$. Hence, $\beta _{pel}$ and $\gamma _{pel}$ are two  quasigeodesics in $ \TC ( X)$ without backtracking joining the same pair of points. By Lemma \ref{pel-track},  $\beta _{pel}$ and $\gamma _{pel}$ have similar intersection patterns with the maximal horospheres-like sets $C_{\alpha}$. Let $\gamma$ be the corresponding electroambient quasigeodesics in $\GXH$, where $\GXH$ is a hyperbolic metric space obtained from $X$ by attaching the hyperbolic cone $C_{\alpha}^h$ to $C_\alpha$, then $\beta$ lies in a bounded neighborhood of $\gamma$.  Due to similar intersection patterns of $\beta _{pel}$ and $\gamma _{pel}$, if $\gamma _{pel}$ penetrates a horosphere $C_\alpha$ and $\beta _{pel}$ does not, then the length of geodesic traversed by $\gamma$ inside $C_{\alpha}$ is uniformly bounded.

Thus there exists $C_1\geq 0$ such that if $x\in \beta ^b_{pel} = \beta _{pel}\setminus \CC$,
then there exists $y\in \gamma ^b_{pel} = \gamma _{pel}\setminus \CC$ such that
$d(x,y)\leq C_1$.

Since $\gamma ^b _{pel}\subset B_{\lambda ^b_{v_0}}$, it follows
from  lemma \ref{final lem}  that \\
$d(y,p)\geq \frac {n}{C+1}$.

So, $\frac {n}{C+1}\leq d(x,y) + d(x,p)\leq C_1 + d(x,p)$,

i.e. $d(x,p)\geq \frac {n}{C+1} - C_1$ (=$M(n)$,say).

Thus we have the following proposition :

\begin{prop}\label{imp thm}
For every point $x$ on $\beta ^b $,  $x$ lies outside an $M(n)$-ball around $p$ in $X$,  such that $M(n)\rightarrow \infty$ as $n\rightarrow \infty$.
\end{prop}

It is now easy to assemble the pieces to deduce the   existence of Cannon-Thurston maps.

\begin{theorem} \label{main thm}
Let $P\colon X\to T$ be a tree of relatively hyperbolic spaces satisfying the quasi-isometrically (qi) embedded condition.
Further suppose that the inclusion of edge-spaces into vertex spaces is strictly type-preserving, and  the induced tree of coned-off spaces continue to satisfy the qi-embedded condition. If $X$ is strongly hyperbolic relative to the family $\CC$ of maximal cone-subtrees of horosphere-like sets, then a Cannon-Thurston map exists for the proper embedding $i_{v_0}\colon X_{v_0}\to X$, where $v_0$ is a vertex of $T$ and $(X_{v_0},d_{X_{v_0}})$ is the relatively hyperbolic metric space corresponding to $v_0$.
\end{theorem}

\noindent {\bf Proof:} A Cannon-Thurston map exists if it satisfies the condition of Lemma \ref{crit-tree}.

So for a fixed reference point $p\in X_{v_0}$ with $p$ lying outside horosphere-like sets,  we assume that $\hat \lambda _{v_0}$ is
 an electric geodesic in $ \hhat X_{v_0}$ such that $\lambda ^b_{v_0} \subset
X_{v_0}$ lies outside an $n$-ball $B_n(p)$ around $p$.

Since $i_{v_0}$ is  a proper embedding,  $\lambda ^b_{v_0}=\hat \lambda_{v_0}\setminus
\bigcup_{H_\alpha\in \HH_{v_0}}H_\alpha$  lies outside an $f(n)$-ball around $p$ in $X$ such that $f(n)\rightarrow \infty$ as $n\rightarrow \infty$.

From the Proposition \ref{imp thm},  if $\beta_{pel}$ is an electrocuted geodesic joining the end points of $\lambda _{v_0}$,  then $\beta ^b$ lies outside an $M(f(n))$-ball around $p$ in $X$ such that $M(f(n))\rightarrow \infty$ as $n\rightarrow \infty$.

From Lemma \ref{crit-tree},
a Cannon-Thurston map for $i\colon X_{v_0}\to X$ exists.
$\Box$

\subsection{Applications and Examples}
In \cite{mitra-trees}, it was shown that for $M$ a hyperbolic 3-manifold of {\bf bounded geometry}
without parabolics, the lifts of
 simply degenerate ends
to the universal cover $\til{M}$ are (uniformly)
quasi-isometric to trees (in fact, rays) of hyperbolic metric spaces satisfying the qi-embedded condition. When $M$
has parabolics, the same arguments show that $\til{M}$, with
cusps excised,  is quasi-isometric to a tree of {\it relatively} hyperbolic metric spaces satisfying the
qi-embedded condition and the conditions of Theorem
\ref{main thm}. Thus, we obtain the following Theorem
of Bowditch:

\begin{theorem} ({\bf Bowditch \cite{bowditch-ct}})
Let $M$ be a 3-manifold corresponding to a Kleinian surface group without accidental parabolics. Further, suppose that
$M$ has bounded geometry. If $S$ denotes the corresponding finite volume hyperbolic surface with some hyperbolic structure, then the inclusion $i: \til{S} \rightarrow
\til{M}$ extends continuously to the boundary, i.e. has a
Cannon-Thurston map.
\label{bowditch}
\end{theorem}

In fact, a more general theorem regarding existence of Cannon-Thurston
maps was shown in \cite{brahma-pared}
for $M$ any bounded
geometry hyperbolic 3 manifold with core incompressible
away from cusps (not necessarily a surface group).
Under the additional assumption that $M$ has no accidental parabolics we conclude that
$\widetilde M$ with cusps excised is quasi-isometric to a tree of relatively
hyperbolic metric spaces. Thus we have:

\begin{theorem} \cite{brahma-pared}
Let $M$ be a hyperbolic
3-manifold corresponding to a Kleinian group of bounded geometry. Further suppose that the compact core
of $M$ is incompressible away from cusps and has no accidental parabolics. If $N$ denotes
a geometrically finite hyperbolic
3 manifold with some hyperbolic structure identified with
the convex core of $M$, then
the inclusion $i: \til{N} \rightarrow
\til{M}$ preserves extends continuously to the boundary, i.e. has a
Cannon-Thurston map.
\label{pared}
\end{theorem}

Much of the strife in \cite{brahma-pared} came from the presence of accidental parabolics
and is in fact the major focal point in Section 4 of that paper. It would be interesting
to generalize the results of this paper to this context.

In \cite{mahan-reeves}, we had, generalizing a Theorem of
Mosher \cite{mosher-hbh}, given relatively hyperbolic examples $G$
of the form

\begin{center}
$1 \rightarrow \pi_1(S) \rightarrow G \rightarrow F_k
\rightarrow 1$
\end{center}

where $F_k$ is free and $S$ is a surface with punctures such that each diffeomorphism corresponding to an element
of $F_k$ preserves the conjugacy class of the puncture. It follows that such
 pairs $(\pi_1(S),G)$ have Cannon-Thurston maps.


\end{document}